\documentclass[12pt]{article}

\usepackage{amsmath, epsfig, cite, lineno}
\usepackage{amssymb,amsthm}
\usepackage{amsfonts, color}
\usepackage{latexsym}

\newtheorem{thm}{Theorem}[section]

\newtheorem{lem}[thm]{Lemma}

\newcommand{\pf}{\noindent{\it Proof.} }

\def\N{{\mathbb N}}

\numberwithin{equation}{section}

\begin{document}


\begin{center}
{\Large\bf On the divisibility of sums involving powers \\[5pt]of multi-variable Schmidt polynomials}
\end{center}

\vskip 2mm \centerline{Qi-Fei Chen$^1$ and Victor J. W. Guo$^2$\footnote{Corresponding author.}}
\begin{center}
{\footnotesize Department of Mathematics, Shanghai Key Laboratory of
PMMP, East China Normal University,\\ 500 Dongchuan Rd., Shanghai
200241,
 People's Republic of China\\
{$^1$\tt figo116@outlook.com,} \quad $^2${\tt jwguo@math.ecnu.edu.cn,\quad
http://math.ecnu.edu.cn/\textasciitilde{jwguo}}}
\end{center}


\vskip 0.7cm \noindent{\bf Abstract.} The multi-variable Schmidt polynomials
are defined by
$$
S_n^{(r)}(x_0,\ldots,x_n):=\sum_{k=0}^n {n+k \choose 2k}^{r}{2k\choose k} x_k.
$$
We prove that, for any positive integers $m$, $n$, $r$, and $\varepsilon=\pm 1$,
all the coefficients in the polynomial
$$
\sum_{k=0}^{n-1}\varepsilon^k(2k+1) S_k^{(r)}(x_0,\ldots,x_k)^m
$$
are multiples of $n$. This generalizes a recent result of Pan on the divisibility of sums of Ap\'ery polynomials.



\section{Introduction}
Schmidt \cite{Sc} introduced the numbers $$\sum_{k=0}^n{n\choose k}^r{n+k\choose k}^r,$$
which, for $r=2$, are called the Ap\'ery numbers since Ap\'ery \cite{Apery} published his ingenious proof of the irrationality
of $\zeta(3)$. Motivated by the work of Z.-W. Sun \cite{Sun1} on the Ap\'ery polynomials,
Guo and Zeng \cite{GZ2} introduced
the Schmidt polynomials
\begin{align*}
S_n^{(r)}(x)=\sum_{k=0}^{n}{n\choose k}^r{n+k\choose k}^r x^k,
\end{align*}
and obtained the following congruence:
$$
\sum_{k=0}^{n-1}\varepsilon^k(2k+1)S_k^{(r)}(x)\equiv0\pmod{n},
$$
where $\varepsilon=\pm 1$. Here and in what follows, we say that a polynomial $P$ (perhaps in multi-variables) is congruent to $0$ modulo $n$,
if all the coefficients in $P$ are multiples of $n$.

Recently, Pan \cite{Pan} proved that, for any positive integers $r$, $m$ and $n$, there holds
\begin{align}
\sum_{k=0}^{n-1}\varepsilon^k(2k+1)S_k^{(r)}(x)^m\equiv0\pmod{n}, \label{eq:pan}
\end{align}
where $\varepsilon=\pm 1$, and thus proved the related conjectures by Sun \cite{Sun1,Sun2} and Guo and Zeng \cite{GZ2}.
Note that, it is not easy to prove \eqref{eq:pan} for $m\geqslant 2$ directly, because there are no
simple formula to calculate the coefficients in the left-hand side of \eqref{eq:pan}.
In fact, Pan \cite{Pan} proved  \eqref{eq:pan} by establishing a $q$-analogue.

Let
$$
S_{n}^{(r)}(x_0,\ldots,x_{n})=\sum_{k=0}^n {n+k \choose 2k}^{r}{2k\choose k} x_k.
$$
We shall give a generalization of \eqref{eq:pan} as follows.

\begin{thm}\label{thm:schmidt}
For any positive integers $r$, $m$, $n$ and nonnegative integer $a$, we have
\begin{align}
\sum_{k=0}^{n-1}(2k+1)S_{k}^{(r)}(x_0,\ldots,x_{k})^m &\equiv0\pmod{n},  \label{eq:main3}\\
\sum_{k=0}^{n-1}(-1)^k (2k+1) S_{k}^{(r)}(x_0,\ldots,x_{k})^m &\equiv0\pmod{n}. \label{eq:main4}
\end{align}
\end{thm}
It is clear that, if we take $x_k={2k\choose k}^{r-1}x^k$ in the congruences \eqref{eq:main3} and \eqref{eq:main4},
then we obtain the congruence \eqref{eq:pan}.

\section{Proof of Theorem \ref{thm:schmidt} for $m=1$}
We first establish the following lemma.
\begin{lem}\label{lem:one} For $\ell,m\in \N$ and $r\geqslant 1$, there exist integers $b_{m,k}^{(r)}$ $(m\leqslant k \leqslant rm)$
divisible by ${k\choose m}$ such that
\begin{align}
{\ell+m\choose 2m}^r{2m\choose m}= \sum_{k=m}^{rm}b_{m,k}^{(r)}{\ell+k\choose 2k}{2k\choose k}.\label{eq:main5}
\end{align}
\end{lem}
\pf We proceed by mathematical induction on $r$. It is clear that the identity \eqref{eq:main5} holds for $r=1$ with $b_{m,m}^{(1)}=1={m\choose m}$.
Suppose that the statement is true for $r$. A special case of the Pfaff-Saalsch\"utz identity
(see  \cite[p.~44, Exercise 2.d]{Stanley} and \cite{Andrews})  reads
$$
{\ell+m\choose m+k}{\ell+k\choose k}= \sum_{i=0}^{k} {\ell-m\choose i}{m\choose m+i-k}{\ell+m+i\choose \ell},
$$
which can be rewritten as
\begin{align}
{\ell+k\choose 2k}{2k\choose k}=\sum_{i=0}^{k} \frac{(m+k)!i!}{(m+i)!k!}{m\choose k-i}{\ell-m\choose i}{\ell+m+i\choose i}.\label{eq:main8}
\end{align}
Multiplying both sides of \eqref{eq:main5} by ${\ell+m\choose 2m}$ and applying \eqref{eq:main8}, we obtain
\begin{align}
&\hskip -2mm {\ell+m\choose 2m}^{r+1}{2m\choose m} \nonumber\\
&= \sum_{k=m}^{rm}b_{m,k}^{(r)}{\ell+k\choose 2k}{2k\choose k}{\ell+m\choose 2m} \nonumber\\
&= \sum_{k=m}^{rm}b_{m,k}^{(r)}\sum_{i=0}^{k} \frac{(m+k)!i!}{(m+i)!k!}{m\choose k-i}
{\ell-m\choose i}{\ell+m+i\choose i}{\ell+m\choose 2m}  \nonumber\\
&= \sum_{k=m}^{rm} \sum_{i=0}^{k} b_{m,k}^{(r)}{k\choose m}^{-1}
{m+i\choose m}{m\choose k-i}{m+k\choose 2m}{\ell+m+i\choose 2(m+i)}{2(m+i)\choose m+i}.  \label{eq:rec-paff}
\end{align}
Letting $m+i=j$ and exchanging the summation order in the right-hand side of \eqref{eq:rec-paff}, we see that
the identity \eqref{eq:main5} holds for $r+1$ with
\begin{align}
b_{m,j}^{(r+1)}={j\choose m}\sum_{k=m}^{rm} b_{m,k}^{(r)}{k\choose m}^{-1}{m\choose j-k}{m+k\choose 2m},
\quad m\leqslant j\leqslant (r+1)m.\label{eq:main9}
\end{align}
By the induction hypothesis, each $b_{m,k}^{(r)}$ is divisible by ${k\choose m}$,
and so $b_{m,j}^{(r+1)}$ is an integer divisible by ${j\choose m}$. This proves that the identity \eqref{eq:main5} is also true for $r+1$.
\qed

\medskip
\noindent{\it Proof of Theorem {\rm\ref{thm:schmidt}} for $m=1$.} Applying \eqref{eq:main5} and the identity
\begin{align}
\sum_{\ell=k}^{n-1}(2\ell+1){\ell+k\choose 2k}{2k\choose k}=n{n\choose k+1}{n+k\choose k},\label{eq:main12}
\end{align}
we have
\begin{align*}
\sum_{k=0}^{n-1}(2k+1) S_{k}^{(r)}(x_0,\ldots,x_{k})
&=\sum_{k=0}^{n-1}(2k+1) \sum_{i=0}^k {k+i\choose 2i}^r{2i\choose i} x_i \\
&=\sum_{i=0}^{n-1} x_i \sum_{k=i}^{n-1} \sum_{j=i}^{ri} b_{i,j}^{(r)}(2k+1){k+j\choose 2j}{2j\choose j}\\
&=n\sum_{i=0}^{n-1} x_i \sum_{j=i}^{ri} b_{i,j}^{(r)}{n\choose j+1}{n+j\choose j}.
\end{align*}
Similarly, applying \eqref{eq:main5} and the identity
\begin{align}
\sum_{\ell=k}^{n-1}(-1)^{\ell}(2\ell+1){\ell+k\choose 2k}{2k\choose k}=(-1)^{n-1}n{n-1\choose k}{n+k\choose k},\label{eq:main13}
\end{align}
we get
\begin{align*}
\sum_{k=0}^{n-1}(-1)^k(2k+1) S_{k}^{(r)}(x_0,\ldots,x_{k})
=(-1)^{n-1}n\sum_{i=0}^{n-1} x_i \sum_{j=i}^{ri} b_{i,j}^{(r)}{n-1\choose j}{n+j\choose j}.
\end{align*}
This completes the proof. \qed

\section{Proof of Theorem \ref{thm:schmidt} for $m>1$}
It is easy to see that
\begin{align}
\sum_{k=0}^{n-1}\varepsilon^k(2k+1)S_{k}^{(r)}(x_0,\ldots,x_{k})^m   &=\sum_{k=0}^{n-1}\varepsilon^k(2k+1)\left( \sum_{i=0}^k {k+i \choose 2i}^{r}{2i\choose i} x_i \right)^m \notag\\
&=\sum_{k=0}^{n-1}\varepsilon^k(2k+1)
\sum_{0\leqslant i_1,\ldots,i_m\leqslant k}
\prod_{j=1}^{m}{k+i_j \choose 2i_j}^{r}{2i_j\choose i_j} x_{i_j}, \label{eq:final}
\end{align}
Exchanging the summation order, we may rewrite the right-hand side of \eqref{eq:final} as
\begin{align}
\sum_{0\leqslant i_1,\ldots,i_m\leqslant n-1}  \sum_{k=\max\{i_1,\ldots,i_m\}}^{n-1}
\varepsilon^k(2k+1) \prod_{j=1}^{m}{k+i_j \choose 2i_j}^{r}{2i_j\choose i_j}x_{i_j}.\label{eq:main11}
\end{align}
By Lemma \ref{lem:one},\ we see that ${k+i_j \choose 2i_j}^{r}{2i_j\choose i_j}$, as a polynomial in $k$, is a linear combination of
$$
{k+i_j \choose 2i_j}{2i_j\choose i_j}, {k+i_j+1\choose 2i_j+2}{2i_j+2\choose i_j+1},\ldots,\
{k+ri_j \choose 2ri_j}{2ri_j\choose ri_j}
$$
with integer coefficients.  Note that \eqref{eq:main8} with $m=j$ and $k=i$ also gives
\begin{align}
{\ell+i\choose 2i}{2i\choose i}{\ell+j\choose 2j}{2j\choose j}
=\sum_{k=0}^{i}{i+j\choose i}{j\choose i-k}{j+k\choose k}{\ell+j+k\choose 2j+2k}{2j+2k\choose j+k}.  \label{eq:repeat}
\end{align}
Therefore, by applying Lemma \ref{lem:one} and then repeatedly using \eqref{eq:repeat}, we find that
the expression $\prod_{j=1}^{m}{k+i_j \choose 2i_j}^{r}{2i_j\choose i_j}$, as a polynomial in $k$, can be written as
a linear combination of
$$
{k+i \choose 2i}{2i\choose i}, {k+i+1\choose 2i+2}{2i+2\choose i+1},\ldots,\
{k+r(i_1+\cdots+i_m) \choose 2r(i_1+\cdots+i_m)}{2r(i_1+\cdots+i_m)\choose r(i_1+\cdots+i_m)},
$$
with integer coefficients, where $i=\max\{i_1,\ldots,i_m\}$.

Finally, using \eqref{eq:main12} and \eqref{eq:main13}, we deduce that
the coefficients in \eqref{eq:main11} are all multiples of $n$. Namely, Theorem \ref{thm:schmidt} holds for $m>1$.

\medskip
\noindent{\it Remark.} The identity \eqref{eq:repeat} was also utilized by the second author
\cite{Guo} to confirm some conjectures of Z.-W. Sun \cite{Sun1,Sun2} on sums of powers of Delannoy polynomials.


\section{Concluding remarks}
Theorem \ref{thm:schmidt} can be further generalized as follows. Let $(x)_0=1$ and $(x)_n=x(x+1)\cdots (x+n-1)$ for all $n\geqslant 1$. Then, for any nonnegative integer $a$,
the proof of \cite[Lemma 2.1]{GZ} gives
\begin{align}
k^a(k+1)^a{k+j\choose 2j}=\sum_{i=0}^{a}c_i(j,a){k+j+i\choose j+i}(2j+1)_{2i},\label{eq:main14}
\end{align}
where $c_0(j,a),\dots,c_a(j,a)$ are integers independent of $k$.

Using the identity \eqref{eq:main14} and the same idea in the previous section, we can show that
(see the proof of \cite[Theorem 1.2]{GZ})
\begin{align*}
\sum_{k=0}^{n-1}\varepsilon^k (2k+1)k^a(k+1)^a S_{k}^{(r)}(x_0,\ldots,x_{k})^m &\equiv0\pmod{n},
\end{align*}
where $\varepsilon=\pm 1$.
Noticing that
\begin{align}
(2k+1)^{2a}=(4k^2+4k+1)^a=\sum_{i=0}^{a}{a\choose i}4^ik^i(k+1)^i,  \label{eq:2a+1}
\end{align}
we immediately get
\begin{align*}
\sum_{k=0}^{n-1}\varepsilon^k(2k+1)^{2a+1}S_{k}^{(r)}(x_0,\ldots,x_{k})\equiv0\pmod{n}, 
\end{align*}
where $\varepsilon=\pm 1$.

\vskip 5mm \noindent{\bf Acknowledgments.} This work was partially
supported by the Fundamental Research Funds for the Central
Universities and the National Natural Science Foundation of China
(grant 11371144).


\begin{thebibliography}{99}
\small \setlength{\itemsep}{-.8mm}

\bibitem{Andrews}G.E. Andrews, Identities in combinatorics, I: on sorting two ordered sets, Discrete Math. 11 (1975), 97--106.

\bibitem{Apery}R. Ap\'ery, Irrationalit\'e de $\zeta(2)$ et $\zeta(3)$, Ast\'erisque 61 (1979), 11--13.

\bibitem{Guo}V.J.W. Guo, Some congruences involving powers of Delannoy polynomials, preprint, 2014.

\bibitem{GZ}V.J.W. Guo and J. Zeng, New congruences for sums involving Apery numbers or central Delannoy numbers,
preprint, Int. J. Number Theory 8 (2012), 2003--2016.

\bibitem{GZ2}V.J.W. Guo and J. Zeng, Proof of some conjectures of Z.-W. Sun on congruences for Ap\'ery polynomials,
J. Number Theory, 132 (2012), 1731--1740.

\bibitem{Pan}Hao Pan, On divisibility of sums of Ap\'ery polynomials, J. Number Theory 143 (2014), 214--223.

\bibitem{Sc}A.L. Schmidt, Generalized $q$-Legendre polynomials,
J. Comput. Appl. Math. 49 (1993), 243--249.

\bibitem{Stanley}R.P. Stanley, Enumerative Combinatorics, Vol. I, Cambridge University Press, Cambridge, 1997.

\bibitem{Sun1}Z.-W. Sun, On sums of Ap\'ery polynomials and related congruences, J. Number Theory, 132 (2012),
2673--2699.

\bibitem{Sun2}Z.-W. Sun, Congruences involving generalized trinomial coefficients, Sci. China 57 (2014), 1375-1400.


\end{thebibliography}
\end{document}